\documentclass[psamsfonts,12pt]{amsart}
\usepackage{amsmath, amssymb, amsthm, latexsym}
\usepackage{hyperref}
\usepackage{amssymb}
\usepackage{latexsym}
\usepackage[center]{caption}
\usepackage{tikz}
\newcounter{braid}
\newcounter{strands}
\pgfkeyssetvalue{/tikz/braid height}{1cm}
\pgfkeyssetvalue{/tikz/braid width}{1cm}
\pgfkeyssetvalue{/tikz/braid start}{(0,0)}
\pgfkeyssetvalue{/tikz/braid colour}{black}
\pgfkeys{/tikz/strands/.code={\setcounter{strands}{#1}}}

\makeatletter
\def\cross{%
  \@ifnextchar^{\message{Got sup}\cross@sup}{\cross@sub}}

\def\cross@sup^#1_#2{\render@cross{#2}{#1}}

\def\cross@sub_#1{\@ifnextchar^{\cross@@sub{#1}}{\render@cross{#1}{1}}}

\def\cross@@sub#1^#2{\render@cross{#1}{#2}}

\def\render@cross#1#2{
  \def\strand{#1}
  \def\crossing{#2}
  \pgfmathsetmacro{\cross@y}{-\value{braid}*\braid@h}
  \pgfmathtruncatemacro{\nextstrand}{#1+1}
  \foreach \thread in {1,...,\value{strands}}
  {
    \pgfmathsetmacro{\strand@x}{\thread * \braid@w}
    \ifnum\thread=\strand
    \pgfmathsetmacro{\over@x}{\strand * \braid@w + .5*(1 - \crossing) * \braid@w}
    \pgfmathsetmacro{\under@x}{\strand * \braid@w + .5*(1 + \crossing) * \braid@w}
    \draw[braid] \pgfkeysvalueof{/tikz/braid start} +(\under@x pt,\cross@y pt) to[out=-90,in=90] +(\over@x pt,\cross@y pt -\braid@h);
    \draw[braid] \pgfkeysvalueof{/tikz/braid start} +(\over@x pt,\cross@y pt) to[out=-90,in=90] +(\under@x pt,\cross@y pt -\braid@h);
    \else
    \ifnum\thread=\nextstrand
    \else
     \draw[braid] \pgfkeysvalueof{/tikz/braid start} ++(\strand@x pt,\cross@y pt) -- ++(0,-\braid@h);
    \fi
   \fi
  }
  \stepcounter{braid}
}

\tikzset{braid/.style={double=\pgfkeysvalueof{/tikz/braid colour},double distance=1pt,line width=2pt,white}}

\newcommand{\braid}[2][]{%
  \begingroup
  \pgfkeys{/tikz/strands=2}
  \tikzset{#1}
  \pgfkeysgetvalue{/tikz/braid width}{\braid@w}
  \pgfkeysgetvalue{/tikz/braid height}{\braid@h}
  \setcounter{braid}{0}
  \let\sigma=\cross
  #2
  \endgroup
}
\makeatother

\input xypic

\swapnumbers
\newtheorem{theorem}[subsection]{Theorem}
\newtheorem{proposition}[subsection]{Proposition}
\newtheorem{lemma}[subsection]{Lemma}
\newtheorem{corollary}[subsection]{Corollary}

\def\Proof{\medskip\noindent{\bf Proof: }}

\def\Z{\mathbb{Z}}

\def\R{\mathbb{R}}

\def\md{\mathcal{D}}

\def\qed{\hfill$\square$\medskip}

\def\Zpk{\mathbb{Z}/p^{k}}
\def\Zpk1{\mathbb{Z}/p^{k-1}}

\newcommand{\rref}[1]{(\ref{#1})}

\newcommand{\fracd}[2]{\frac{\displaystyle #1}{\displaystyle #2}}

\newcommand{\beg}[2]{\begin{equation}\label{#1}#2\end{equation}}
\def\r{\rightarrow}
\def\mc{\mathcal{C}}

\def\sl2{\widetilde{SL_{2}(\Z)}}

\def\md{\mathcal{D}}

\usepackage[all,arc]{xy}
\usepackage{enumerate}
\usepackage{mathrsfs}
\title[Some HF and BOS homologies of branched double covers]{Some non-trivial
examples of the Baldwin-Ozsv\'{a}th-Szab\'{o}
twisted spectral sequence and Heegaard-Floer homology of branched double covers}
\author{Elden Elmanto and Igor Kriz}
\thanks{Elden Elmanto was supported in part by the University of Chicago. 
Igor Kriz was supported by NSF grant DMS 1102614, and the Wolfensohn fund at the 
Institute for Advanced Study}
\begin{document}

\begin{abstract}
We present some non-trivial calculations of Baldwin-Ozsv\'{a}th-Szab\'{o} cohomology of links, and applications to Heegaard-Floer homology of branched double covers.
\end{abstract}

\maketitle

\section{Introduction}

\vspace{3mm}

The main point of the present project was to compute  some non-trivial examples
of Baldwin-Ozsv\'{a}th-Szab\'{o} (BOS) cohomology of links \cite{os,kk}, and,
to apply this to Heegaard-Floer homology of branched double covers
via the Baldwin-Ozsv\'{a}th-Szab\'{o} twisted spectral sequence \cite{os}.
The examples we consider here are of a different type then what the main focus
of recent interest has been, for reasons we shall explain.

Much recent interest has focused on so called L-spaces, which are
compact oriented $3$-manifolds satisfying
\beg{elsp}{rank(\widehat{HF}(Y)=|H_1(Y,\Z)|.}
For example, Lisca and Stipsicz \cite{ls1,ls2} characterized L-spaces which are
Seifert-fibered over $S^2$. For a branched double cover $\Sigma(L)$ of a link $L$,
condition \rref{elsp} is equivalent to
\beg{elsp1}{rank(\widehat{HF}(\Sigma(L))=det(L)\neq 0.
} 
Baldwin, Ozsv\'{a}th and Szab\'{o} \cite{os} constructed a twisted variant
of their spectral sequence \cite{os2} convergent to $\widehat{HF}(\Sigma(L))$,
which we will describe below. The point is that the (combinatorially
defined) $E_3$-term of the twisted
spectral sequence, which we call Baldwin-Ozsv\'{a}th-Szab\'{o} (briefly BOS)
cohomology, is extremely sparse (to the point that one may conjecture
it collapses, although that is not known at present). In \cite{kk}, it was proved that BOS cohomology
is an  invariant of oriented links, and it was observed that this gives an good method for
detecting links which satisfy \rref{elsp1}. The reason is that while BOS cohomology
is defined as the cohomology of a cochain complex defined over a field of rational
functions over $\mathbb{F}_2$ in a number of variables which increases with 
the number of crossings of $L$ (described below), from an algebraic point of
view, \rref{elsp1} is ``generic behavior'' on such complexes. Roughly,
``generic behavior'' menas that anything that can cancel cancels.
For example in a $2$-stage complex (given by a single
linear map), generic behavior means that the rank of the linear map is maximal
allowed by the dimension. As explained in \cite{kk} (and also used in
this note), there is a mehod for detecting generic behavior in BOS cohomology.
Namely, it is possible to set all the variables
in the fraction field equal to integral powers of a single variable. If one gets lucky
and the $E_3$-term has rank equal to $det(L)$, it is also true in BOS cohomology,
and in $\widehat{HF}$.
While calculations in fields of rational functions in many variables are computationally
extremely inefficient, calculations in one variable are no problem. Because of this,
\rref{elsp1} can be detected by BOS cohomology, and this was used
in \cite{kk} to find a new weaker condition on links whose branched double covers
are L-spaces.

\vspace{3mm}
To complement this, in the present project, we wanted to 
do computations where BOS cohomology behaves
non-generically, with some applications to Heegaard-Floer
homology. Since there is an algorithm \cite{osp} for calculating 
Heegaard-Floer homology of Seifert-fibered spaces,
examples of hyperbolic knots are of most interest.
We found that despite the combinatorial definition, 
it is {\em extraordinarily difficult} to compute directly non-generically behaved
examples of BOS cohomology. This note serves, perhaps, as
a case study of the difficulty of such computations.
In the end, combining heuristics with computer-assisted
methods, we succeeded in computing one example (perhaps one of
the smallest ones) in Proposition \ref{p1} below. The example
happens to be a link with $0$ determinant, and infinitely many
examples of exact computation of BOS cohomology can be deduced 
using the skein behavior of BOS cohomology (Proposition \ref{p1s2} and Corollary \ref{c1}).
Since all these examples are links of determinant $0$, we do not get
an immediate application to $\widehat{HF}$. However, using the 
Ozsv\'{a}th-Szab\'{o} computation of $\widehat{HF}$ of
$T(7,3)$ (\cite{ost3}) as input,
we were able to calculate BOS-cohomology and $\widehat{HF}$ for
infinitely many new examples (Theorem \ref{t1}, Theorem \ref{t4}),
all but finitely many of which are hyperbolic (Proposition \ref{thyp}).
The reader should keep in mind that although we have no example
of non-collapse of the twisted BOS spectral sequence, BOS cohomology
carries more information than $\widehat{HF}$, since
non-trivial ranks appear in different degrees; because of this, 
the BOS calculations are also of independent interest. 
We consider our examples as a ``proof of concept'',
showing how this method works; it is very likely that many other examples
can be calculated in a similar way. At the same time, it is also clear
that such examples do not come cheap.

\vspace{3mm}
The present paper is organized as follows. In Section \ref{sp} we review the
preliminaries, i.e. the definition of BOS cohomology, the link invariance of
\cite{kk} and the BOS spectral sequence to $\widehat{HF}$. In Section
\ref{s1}, we treat the one example of non-trivial BOS cohomology which we
were able to compute directly. In Section \ref{s2}, we treat all the examples
derived from this and from what was known about $T(7,3)$ in \cite{ost3}.

\vspace{3mm}

\noindent
{\bf Acknowledgement:} We are very indebted to P.Ozsv\'{a}th for comments.

\vspace{3mm}

\section{Preliminaries}
\label{sp}

\vspace{3mm}
Let us first describe BOS cohomology. Let $\md=\md(L)$ be a non-degenerate 
projection of a link $L$, i.e. an embedding of a link in $\R^3$ such that the projection
on the $xy$ plane is an immersion with at most finitely many double crossings.
Then the faces of the projection can be colored black and white so
that no two faces of the same color border the same arc of the projection.
Form a planar graph (with possible multiple edges and loops)
whose vertices are the faces colored black, and edges are crossings
which border two faces colored black. This is called the 
{\em black graph} $B(\md)$. An edge has {\em height $0$}
if the  edge of the black graph, the lower arc of the crossing and the upper
arc of the crossing occur in this order clockwise, and {\em height $1$} otherwise.

\vspace{3mm}
The BOS cochain complex $\mc_0(\md)$  is then formed as follows: 
Pick one vertex of the black graph as a base point. Let $F$ be
a field of rational functions over $\mathbb{F}_2$ on variables
corresponding to all vertices of the black graph other than the base
point and all bounded faces of the black graph.
$\mc_0(\md)_k$ is the free $F$-module on all spaning trees $T$ of $B(\md)$
whose total height is $2k$; here the total height is defined to be the
number of edges of height $1$ included in $T$ plus the number
of edges of height $0$ not included in $T$.
The number $k$ can be a half-integer, but all the possible values of $k$
differ by an integer.

\vspace{3mm}
The differential $\Psi$ of $\mc_0$ increases total degree by $1$. A non-trivial
coefficient occurs between a spanning tree $T$ and a spanning tree $T^\prime$
obtained from $T$ by removing one edge $e$ of height $0$ and adding one edge $f$ of
height $1$.
The coefficient is of the form
$$\frac{1}{1+\alpha} + \frac{1}{1+\beta},$$
where $\alpha$ is the product of all the variables corresponding
to faces enclosed inside the cycle $c$ in $T\cup \{f\}$
provided that $e$, $f$ and the base point occur counter-clockwise on that cycle, and the product
of the inverses of all the variables corresponding to faces enclosed inside of $c$ otherwise;
the element $\beta$ is the product of all the variables corresponding to vertices 
in the component of $T\smallsetminus \{e\}$ not containing the base point.

It was proved directly in \cite{kk} that $\Psi$ is a differential,
i.e. that $\Psi\circ\Psi=0$. A key result of \cite{os} is the following

\vspace{3mm}
\begin{theorem}
\label{tos}
Let $L$ be a link with non-zero determinant. Then there exists a single-graded spectral
sequence
$$E_3=H^*(\mc_0(\md(L)))\Rightarrow \widehat{HF}(\Sigma(L))\otimes_{\mathbb{F}_2}F.$$
Moreover, the grading is by total height (i.e. twice the degree), and the spectral
sequence is sparse in the sense that the only possible non-zero differentials
are of the form $d_{4k+2}$, $k\in\Z$.
\end{theorem}

\vspace{3mm}
To make BOS cohomology a knot invariant, one must correct by half the number of
negative crossings, to take care of Reidemeister 1 moves. For an oriented link, a postive
crossing is one where the upper arc of the crossing goes from lower left to upper right
and the lower arc goes from lower right to upper left. The other kind of crossing is
called a negative crossing. Let
$$\mc(\md)_k=\mc_0(\md)_{k+n_-/2}$$
where $n_-$ is the number of negative crossings in $\md$. The differential in
$\mc$ is defined to be the same as in $\mc_0$. 
In \cite{kk}, the following was proved:

\vspace{3mm}
\begin{theorem}
\label{tkk}
The numbers
$$rank H^i(\mc(\md))$$
are invariants of oriented links and unoriented knots.
\end{theorem}

We therefore put
$$H^{i}_{BOS}(L)=H^i(\mc(\md(L)).$$

\vspace{3mm}
In this paper, we will work with the unshifted BOS cohomology, i.e. the cohomology
of the complex $\mc_0(\md)$ of a projection $\md$ where generating trees are graded
by $1/2$ times the total height. This is only a link invariant up to shift, but if $\md^i$
are the $i$-resolutions of $\md$ along a single edge, $i=0,1$, we have a long
exact sequence
\beg{ei1}{\dots\r H^{i-1}(\mc_0(\md^0))\r H^{i-1/2}(\mc_0(\md^1))
\r H^i(\mc_0(\md))\r H^i(\mc_0(\md^0))\r \dots
}

\vspace{3mm}
Denote by $B_{n_1,\dots,n_k}=B(\md_{(n_1,\dots,n_k)})$ the 
black graph in Figure 1.

\begin{figure}
\definecolor{c020000}{RGB}{2,0,0}
\definecolor{cffffff}{RGB}{255,255,255}

\begin{tikzpicture}[y=0.60pt, x=0.6pt,yscale=-1, inner sep=0pt, outer sep=0pt]
\path[shift={(38.624983,55.850641)},draw=black,fill=black,line join=miter,line
  cap=butt,even odd rule,line width=0.800pt]
  (340.9265,461.1704)arc(0.000:180.000:6.818530 and
  6.313)arc(-180.000:0.000:6.818530 and 6.313) -- cycle;
\path[draw=black,dash pattern=on 0.54pt off 6.54pt,line join=miter,line
  cap=butt,miter limit=4.00,line width=0.545pt] (221.4478,442.2542) --
  (184.0638,563.1429);
\path[draw=black,dash pattern=on 0.51pt off 6.08pt,line join=miter,line
  cap=butt,miter limit=4.00,line width=0.507pt] (534.4504,441.7410) --
  (570.5268,544.7119);
\path[draw=black,dash pattern=on 7.51pt off 7.51pt,line join=miter,line
  cap=butt,miter limit=4.00,line width=0.626pt] (366.9165,327.6332) .. controls
  (377.8039,321.9039) and (323.4074,348.1617) .. (304.1859,358.7888);
\path[draw=black,dash pattern=on 6.67pt off 6.67pt,line join=miter,line
  cap=butt,miter limit=4.00,line width=0.556pt] (368.7962,323.1577) .. controls
  (415.9620,342.6157) and (447.2796,362.1085) .. (454.4134,360.5334);
\path[draw=black,line join=miter,line cap=butt,line width=1.645pt]
  (370.6080,328.6710) .. controls (370.6080,391.3511) and (370.6080,454.0313) ..
  (370.6080,516.7114);
\path[xscale=0.670,yscale=1.494,draw=black,fill=black,line width=0.389pt]
  (705.4375,244.72665) node[above right] (text3621) {$a_{km_k-1}$     };
\path[xscale=0.690,yscale=1.449,draw=black,fill=black,line width=0.668pt]
  (476.26028,284.53372) node[above right] (text3625) {$e_0$     };
\path[shift={(38.624983,54.840491)},draw=black,fill=c020000,line join=miter,line
  cap=butt,even odd rule,line width=0.800pt]
  (270.7143,302.7193)arc(0.000:180.000:5.535715 and
  5.000)arc(-180.000:0.000:5.535715 and 5.000) -- cycle;
\path[shift={(38.624983,59.891251)},draw=black,fill=black,line join=miter,line
  cap=butt,even odd rule,line width=0.800pt]
  (337.1429,266.1122)arc(0.000:180.000:5.535715 and
  4.821)arc(-180.000:0.000:5.535715 and 4.821) -- cycle;
\path[shift={(38.624983,59.891251)},draw=black,fill=black,line join=miter,line
  cap=butt,fill opacity=0.987,even odd rule,line width=0.800pt]
  (416.0714,301.6479)arc(0.000:180.000:5.892860 and
  5.714)arc(-180.000:0.000:5.892860 and 5.714) -- cycle;
\path[shift={(35.053483,63.462681)},draw=black,fill=black,dash pattern=on 0.80pt
  off 6.40pt,miter limit=4.00,line width=0.800pt]
  (212.1429,330.9336)arc(0.000:180.000:6.071430 and
  5.714)arc(-180.000:0.000:6.071430 and 5.714) -- cycle;
\path[shift={(49.339183,60.605541)},draw=black,fill=black,dash pattern=on 0.80pt
  off 6.40pt,miter limit=4.00,line width=0.800pt]
  (463.5714,327.3622)arc(0.000:180.000:7.142860 and
  6.429)arc(-180.000:0.000:7.142860 and 6.429) -- cycle;
\path[draw=black,dash pattern=on 9.60pt off 9.60pt,line join=miter,line
  cap=butt,miter limit=4.00,line width=0.800pt] (302.9107,357.9678) --
  (238.6249,395.1106);
\path[draw=black,dash pattern=on 9.60pt off 9.60pt,line join=miter,line
  cap=butt,miter limit=4.00,line width=0.800pt] (447.1964,362.2535) --
  (505.7678,389.3963);
\path[xscale=0.828,yscale=1.208,fill=black] (351.09113,274.15787) node[above
  right] (text3643) {$a_{11}$     };
\path[xscale=0.784,yscale=1.276,fill=black] (511.58441,245.57683) node[above
  right] (text3647) {$a_{km_k}$     };
\path[xscale=0.718,yscale=1.393,fill=black] (342.88501,411.01691) node[above
  right] (text3651) {$e_1$     };
\path[shift={(-159.23222,131.56493)},draw=black,fill=black,line join=miter,line
  cap=butt,even odd rule,line width=0.800pt]
  (340.9265,461.1704)arc(0.000:180.000:6.818530 and
  6.313)arc(-180.000:0.000:6.818530 and 6.313) -- cycle;
\path[shift={(250.05348,122.27922)},draw=black,fill=black,line join=miter,line
  cap=butt,even odd rule,line width=0.800pt]
  (340.9265,461.1704)arc(0.000:180.000:6.818530 and
  6.313)arc(-180.000:0.000:6.818530 and 6.313) -- cycle;
\path[draw=black,dash pattern=on 0.74pt off 8.92pt,line join=miter,line
  cap=butt,miter limit=4.00,line width=0.744pt] (204.8610,750.4505) .. controls
  (314.3343,913.5200) and (461.0605,908.2294) .. (559.0150,742.5965);
\path[draw=black,dash pattern=on 9.60pt off 9.60pt,line join=miter,line
  cap=butt,miter limit=4.00,line width=0.800pt] (177.5910,595.2614) --
  (177.5910,670.9757);
\path[draw=black,dash pattern=on 9.60pt off 9.60pt,line join=miter,line
  cap=butt,miter limit=4.00,line width=0.800pt] (583.9688,586.6900) --
  (583.9688,653.8328);
\path[shift={(-158.22202,209.34668)},draw=black,fill=black,line join=miter,line
  cap=butt,even odd rule,line width=0.800pt]
  (340.9265,461.1704)arc(0.000:180.000:6.818530 and
  6.313)arc(-180.000:0.000:6.818530 and 6.313) -- cycle;
\path[shift={(249.87958,196.2147)},draw=black,fill=black,line join=miter,line
  cap=butt,even odd rule,line width=0.800pt]
  (340.9265,461.1704)arc(0.000:180.000:6.818530 and
  6.313)arc(-180.000:0.000:6.818530 and 6.313) -- cycle;
\path[draw=black,dash pattern=on 9.60pt off 9.60pt,line join=miter,line
  cap=butt,miter limit=4.00,line width=0.800pt] (302.9107,357.9678) --
  (238.6249,395.1106);
\path[draw=black,dash pattern=on 9.59pt off 9.59pt,line join=miter,line
  cap=butt,miter limit=4.00,line width=0.799pt] (176.0023,671.3661) --
  (208.3280,760.2603);
\path[draw=black,dash pattern=on 9.26pt off 9.26pt,line join=miter,line
  cap=butt,miter limit=4.00,line width=0.772pt] (587.1523,660.2373) --
  (558.8330,744.1150);
\path[draw=black,line join=miter,line cap=butt,line width=0.800pt]
  (183.0738,666.3157) -- (367.9317,520.8538);
\path[draw=black,line join=miter,line cap=butt,line width=0.800pt]
  (582.0840,657.2243) -- (376.0129,520.8538) -- (525.5155,785.5137);
\path[draw=black,line join=miter,line cap=butt,line width=0.800pt]
  (372.9824,521.8639) -- (480.0586,830.9706);
\path[draw=black,dash pattern=on 0.89pt off 10.68pt,line join=miter,line
  cap=butt,miter limit=4.00,line width=0.890pt] (260.0957,635.6698) .. controls
  (284.5759,768.5951) and (393.5454,752.7871) .. (427.1491,711.5669);
\path[draw=black,dash pattern=on 9.60pt off 9.60pt,line join=miter,line
  cap=butt,miter limit=4.00,line width=0.800pt] (243.6829,394.5847) --
  (225.5002,434.9908);
\path[draw=black,dash pattern=on 9.60pt off 9.60pt,line join=miter,line
  cap=butt,miter limit=4.00,line width=0.800pt] (510.3632,393.5745) --
  (529.5561,430.9502);
\path[draw=black,dash pattern=on 9.60pt off 9.60pt,line join=miter,line
  cap=butt,miter limit=4.00,line width=0.800pt] (176.0027,592.5746) --
  (189.1347,548.1279);
\path[draw=black,dash pattern=on 9.60pt off 9.60pt,line join=miter,line
  cap=butt,miter limit=4.00,line width=0.800pt] (584.1043,580.4528) --
  (568.9520,540.0467);
\path[xscale=0.605,yscale=1.654,fill=black] (376.58566,223.15311) node[above
  right] (text4481) {$a_{12}$     };
\path[xscale=0.609,yscale=1.642,fill=black] (189.2068,389.45865) node[above
  right] (text4485) {$a_{1m_1}$     };
\path[xscale=0.609,yscale=1.642,fill=black] (975.90668,382.28885) node[above
  right] (text4489) {$a_{k1}$     };
\path[shift={(-142.05962,252.78324)},draw=black,fill=black,line join=miter,line
  cap=butt,even odd rule,line width=0.800pt]
  (340.9265,461.1704)arc(0.000:180.000:6.818530 and
  6.313)arc(-180.000:0.000:6.818530 and 6.313) -- cycle;
\path[xscale=0.707,yscale=1.415,fill=black] (167.69576,494.76404) node[above
  right] (text4495) {$a_{21}$     };
\path[xscale=0.718,yscale=1.393,fill=black] (676.27588,413.91846) node[above
  right] (text4499) {$e_{k-1}$     };
\path[shift={(236.74758,250.76293)},draw=black,fill=black,line join=miter,line
  cap=butt,even odd rule,line width=0.800pt]
  (340.9265,461.1704)arc(0.000:180.000:6.818530 and
  6.313)arc(-180.000:0.000:6.818530 and 6.313) -- cycle;
\path[xscale=0.417,yscale=2.397,fill=black] (1384.9164,291.96774) node[above
  right] (text4505) {$a_{k-1m_{k-1}}$     };
\path[draw=black,line join=miter,line cap=butt,line width=0.800pt]
  (408.3378,368.3207) -- (400.2566,337.0060) -- (428.5408,323.8740);
\path[draw=black,line join=miter,line cap=butt,line width=0.800pt]
  (262.8758,357.2090) -- (255.8047,384.4832) -- (279.0383,389.5339);
\path[draw=black,line join=miter,line cap=butt,line width=0.800pt]
  (194.1854,773.3919) -- (224.4900,777.4325) -- (224.4900,749.1482);
\path[draw=black,line join=miter,line cap=butt,line width=0.800pt]
  (238.6322,594.5949) -- (270.9570,595.6050) -- (270.9570,627.9299);
\path[draw=black,line join=miter,line cap=butt,line width=0.800pt]
  (354.7997,450.1431) -- (369.9520,475.3969) -- (390.1550,444.0822);
\path[draw=black,line join=miter,line cap=butt,line width=0.800pt]
  (454.8048,593.5847) -- (451.7743,570.3512) -- (472.9875,563.2802);
\path[draw=black,line join=miter,line cap=butt,line width=0.800pt]
  (409.7339,691.6900) -- (419.0196,657.4043) -- (438.3053,670.9757);
\path[draw=black,line join=miter,line cap=butt,line width=0.800pt]
  (473.3053,726.6900) -- (475.4482,698.8328) -- (501.8767,703.8328);
\path[shift={(1044.7339,-121.38652)},draw=black,fill=cffffff,miter
  limit=4.00,fill opacity=0.000,line width=0.800pt]
  (-645.0000,637.3622)arc(-0.000:180.000:26.785715 and
  22.143)arc(-180.000:0.000:26.785715 and 22.143) -- cycle;

\end{tikzpicture}
\caption{}
\end{figure}

Denote the corresponding link by $L_{(n_1,\dots,n_k)}$. 

\vspace{3mm}

\section{The link $L_{(3,3,0)}$}
\label{s1}

\vspace{3mm}

Our first result is the following

\vspace{3mm}

\begin{proposition}
\label{p1}
We have
$$
rank H^i(\mc_0(\md_{(3,3,0)})=\left\{
\begin{array}{ll}
1 & \text{if $i=1,2$}\\
0 & \text{otherwise.}
\end{array}
\right.
$$

\end{proposition}

\vspace{3mm}

The proof will occupy the
remainder of this section. First note that $B_{(3,3,0)}$ has $18$
spanning trees in heights $2$ and $4$ (and none others). The differential
is, then, an $18\times 18$ matrix $N$ over a field of rational
functions over $\mathbb{F}_2$ with $\geq 6$ variables (it can be reduced from $9$ to $6$ by
the Fundamental lemma of \cite{kk}). It follows then that
$$
rank H^i(\mc_0(\md_{(3,3,0)}))=\left\{
\begin{array}{ll}
r & \text{if $i=1,2$}\\
0 & \text{otherwise}
\end{array}
\right.
$$
for some number $r=0,1,2,\dots$, and
our statement is equivalent
to saying that $r=1$, i.e. that $N$ has rank $17$. This could, in principle,
be checked by computer, but it exceeds the computing power of implementations
of computer algebra softwares we could find. 

\vspace{3mm}

Because of this, we simplified the problem as follows. Consider, instead, $\md_{(3,3)}$
(a projection of the knot $8_{19}$ in Rolfsen's table), and label its vertices and faces
as in Figure 2.

\begin{figure}
\definecolor{cffffff}{RGB}{255,255,255}

\begin{tikzpicture}[y=0.80pt, x=0.8pt,yscale=-1, inner sep=0pt, outer sep=0pt]
\path[draw=black,dash pattern=on 7.51pt off 7.51pt,line join=miter,line
  cap=butt,miter limit=4.00,line width=0.626pt] (1597.0564,784.6877) .. controls
  (1607.9438,778.9584) and (1553.5472,805.2162) .. (1534.3257,815.8433);
\path[draw=black,dash pattern=on 6.67pt off 6.67pt,line join=miter,line
  cap=butt,miter limit=4.00,line width=0.556pt] (1598.9360,780.2122) .. controls
  (1646.1018,799.6702) and (1725.2766,844.8773) .. (1732.4103,843.3022);
\path[xscale=0.821,yscale=1.217,draw=black,fill=black,line width=0.322pt]
  (1825.6709,670.67676) node[above right] (text8750) {$a_{11}$     };

\path[xscale=0.821,yscale=1.217,draw=black,fill=black,line width=0.322pt]
  (1825.6709,934.67676) node[above right] (text8750) {$a_{13}$     };
\path[xscale=0.619,yscale=1.614,draw=black,fill=black,line width=0.290pt]
  (2696.9197,704.04602) node[above right] (text8754) {$a_{21}$     };
\path[shift={(1208.5048,475.027)},draw=black,fill=black,line join=miter,line
  cap=butt,even odd rule,line width=0.800pt]
  (274.2564,582.1362)arc(-0.051:180.051:7.071068 and
  6.566)arc(-180.051:0.051:7.071068 and 6.566) -- cycle;
\path[shift={(1317.6321,463.79277)},draw=black,fill=black,line join=miter,line
  cap=butt,even odd rule,line width=0.800pt]
  (424.2641,580.3684)arc(0.000:180.000:7.828685 and
  6.313)arc(-180.000:0.000:7.828685 and 6.313) -- cycle;
\path[shift={(1268.7648,516.94578)},draw=black,fill=black,line join=miter,line
  cap=butt,even odd rule,line width=0.800pt]
  (337.1429,266.1122)arc(0.000:180.000:5.535715 and
  4.821)arc(-180.000:0.000:5.535715 and 4.821) -- cycle;
\path[shift={(1265.9076,494.08864)},draw=black,fill=black,line join=miter,line
  cap=butt,even odd rule,line width=0.800pt]
  (340.9265,461.1704)arc(0.000:180.000:6.818530 and
  6.313)arc(-180.000:0.000:6.818530 and 6.313) -- cycle;
\path[shift={(1265.1934,520.51721)},draw=black,fill=black,dash pattern=on 0.80pt
  off 6.40pt,miter limit=4.00,line width=0.800pt]
  (212.1429,330.9336)arc(0.000:180.000:6.071430 and
  5.714)arc(-180.000:0.000:6.071430 and 5.714) -- cycle;
\path[shift={(1279.4791,517.66007)},draw=black,fill=black,dash pattern=on 0.80pt
  off 6.40pt,miter limit=4.00,line width=0.800pt]
  (463.5714,327.3622)arc(0.000:180.000:7.142860 and
  6.429)arc(-180.000:0.000:7.142860 and 6.429) -- cycle;
\path[draw=black,dash pattern=on 9.60pt off 9.60pt,line join=miter,line
  cap=butt,miter limit=4.00,line width=0.800pt] (1533.0505,815.0222) --
  (1468.7648,852.1651);
\path[xscale=0.743,yscale=1.346,fill=black] (1902.7317,705.83612) node[above
  right] (text5034) {$a_{12}$     };
\path[xscale=0.786,yscale=1.272,fill=black] (2148.2739,632.07977) node[above
  right] (text5038) {$a_{23}$     };
\path[draw=black,dash pattern=on 9.60pt off 9.60pt,line join=miter,line
  cap=butt,miter limit=4.00,line width=0.800pt] (1473.0505,851.4508) --
  (1475.1934,1052.8794);
\path[draw=black,dash pattern=on 9.60pt off 9.60pt,line join=miter,line
  cap=butt,miter limit=4.00,line width=0.800pt] (1735.1934,846.4508) --
  (1734.4791,1038.5937);
\path[draw=black,dash pattern=on 9.14pt off 9.14pt,line join=miter,line
  cap=butt,miter limit=4.00,line width=0.762pt] (1477.7053,1059.5163) --
  (1603.7188,1167.6075) -- (1734.4269,1046.8453);
\path[draw=black,line join=miter,line cap=butt,line width=0.801pt]
  (1599.0852,783.4098) -- (1601.1044,1165.2464);
\path[shift={(1333.862,582.527)},draw=black,fill=black,line join=miter,line
  cap=butt,even odd rule,line width=0.800pt]
  (274.2564,582.1362)arc(-0.051:180.051:7.071068 and
  6.566)arc(-180.051:0.051:7.071068 and 6.566) -- cycle;
\path[xscale=0.726,yscale=1.377,fill=black] (2397.3735,684.88318) node[above
  right] (text4551) {$a_{22}$     };
\path[xscale=0.877,yscale=1.140,fill=black] (1806.9175,673.33289) node[above
  right] (text4555) {$x$     };
\path[xscale=0.923,yscale=1.084,fill=black] (1715.879,1104.8436) node[above
  right] (text4559) {$v$     };
\path[xscale=0.743,yscale=1.346,fill=black] (1922.1594,631.15784) node[above
  right] (text4563) {$y_1$     };
\path[xscale=0.869,yscale=1.151,fill=black] (1640.2737,925.53436) node[above
  right] (text4567) {$y_2$     };
\path[xscale=0.629,yscale=1.589,fill=black] (2781.0618,530.10394) node[above
  right] (text4571) {$z_1$     };
\path[xscale=0.877,yscale=1.141,fill=black] (1990.41,925.11035) node[above
  right] (text4575) {$z_2$     };
\path[xscale=1.363,yscale=0.733,fill=black] (1104.65,1309.3359) node[above
  right] (text4579) {$f_1$     };
\path[xscale=1.048,yscale=0.954,fill=black] (1564.1765,1005.8182) node[above
  right] (text4583) {$f_2$     };
\path[xscale=0.679,yscale=1.472,fill=black] (2358.4092,599.91028) node[above
  right] (text4587) {$e_0$     };
\path[xscale=0.703,yscale=1.422,fill=black] (2283.8232,748.84259) node[above
  right] (text4591) {$e_1$     };
\path[cm={{0.9478659,0.0,0.0,0.85173729,(1284.6894,577.79635)}},draw=black,fill=cffffff,miter
  limit=4.00,fill opacity=0.000,line width=0.800pt]
  (359.2857,444.1479)arc(-0.022:180.022:26.428572 and
  18.929)arc(-180.022:0.022:26.428572 and 18.929) -- cycle;

\end{tikzpicture}
\caption{}
\end{figure}

Put
$$F_0=\mathbb{F}_2(x,y_1,y_2,z_1,z_2,v),$$
$$F=F_0(T),$$
$$K=F(Q).$$
(The notation means adjoining algebraically independent variables,
i.e. fields of rational functions.)
The variables of $F_0$ are simply the vertex variables which occur in
the definition of BOS cohomology; the variables $T,Q$ are related to the
face variables by
$$f_1f_2=T,\; f_1=Q.$$
It will also be convenient for our purposes to put
$$A=\fracd{1}{1+f_1}=\fracd{1}{1+Q},\;\;
B=\fracd{1}{1+f_2}=\fracd{Q}{Q+T}.$$
Denote by $T_{\epsilon, i,j}$, $\epsilon\in\{0,1\}$, $i\in \{1,2\}$, $j\in \{1,2,3\}$
the spanning tree of $B_{(3,3)}$ obtained by omitting the edge $e_\epsilon$ and
$a_{ij}$. Denote by $T^{\prime}_{i,j}$, $i,j\in\{1,2,3\}$ the spannig tree obtained 
by omitting the edges $a_{1i}$ and $a_{2j}$. 

It will also be convenient to have matrix rows and columns  indexed by simple numbers,
so put
$$\begin{array}{lll}
u_0=T_{0,2,3}, & u_1=T_{0,2,2}, & u_2=T_{0,2,1},\\
u_3=T_{0,1,1}, & u_4=T_{0,1,2}, & u_5=T_{0,1,3},\\
u_6=T_{1,2,3}, & u_7=T_{1,2,2}, & u_8=T_{1,2,1},\\
u_9=T_{1,1,1}, & u_{10}=T_{1,1,2}, & u_{11}=T_{1,1,3}
\end{array}
$$
(these will correspond to columns)
and
$$\begin{array}{lll}
v_1=T^{\prime}_{1,1}, & v_2=T^{\prime}_{2,1}, & v_3=T^{\prime}_{3,1}, \\
v_4=T^{\prime}_{1,2}, &v_5=T^{\prime}_{2,2}, &v_6=T^{\prime}_{3.2}, \\
v_7=T^{\prime}_{1,3}, &v_8=T^{\prime}_{2,3}, &v_9=T^{\prime}_{3,3}
\end{array}$$
(these will correspond to rows) of the $\{1,\dots,9\}\times\{0,\dots,11\}$ matrix $M$
of the differential $\Psi$ of $\mc_0(\md_{(3,3)})$. The entries of the matrix $M$ are given
explicitly by
$$\begin{array}{lll}
M_{1,0}=A+\fracd{1}{1+x}, & M_{2,0}=A+\fracd{1}{1+xy_1}, \\[3ex]M_{3,0}=A+\fracd{1}{1+xy_1y_2},& \\[3ex]
M_{4,1}=A+\fracd{1}{1+xz_1}, & M_{5,1}=A+\fracd{1}{1+xy_1z_1}, \\[3ex]
M_{6,1}=A+\fracd{1}{1+xy_1y_2z_1}, &
\\[3ex]
M_{7,2}=A+\fracd{1}{1+xz_1z_2}, & M_{8,2}=A+\fracd{1}{1+xy_1z_1z_2}, \\[3ex]
M_{9,2}=A+\fracd{1}{1+xy_1p_2z_1z_2},& \end{array}$$
$$\begin{array}{lll}
M_{1,3}=B+1+\fracd{1}{1+x}, & M_{4,3}=B+1+\fracd{1}{1+xz_1}, \\[3ex]
M_{7,3}=B+1+\fracd{1}{1+xz_1z_2}, &\\[3ex]
M_{2,4}=B+1+\fracd{1}{1+xy_1}, & M_{5,4}=B+1+\fracd{1}{1+xy_1z_1}, \\[3ex]
M_{8,4}=B+1+\fracd{1}{1+xy_1z_1z_2}, & \\[3ex]
M_{3,5}=B+1+\fracd{1}{1+xy_1y_2}, & M_{6,5}=B+1+\fracd{1}{1+xy_1y_2z_1}, \\[3ex]
M_{9,5}=B+1+\fracd{1}{1+xy_1y_2z_1z_2}. &
\end{array}
$$
Additionally, the entry $M_{i,6+j}$ is obtained from the entry $M_{i,j}$, $i=0,\dots,5$
by replacing $A$ by $A+1$, $B+1$ by $B$ (to account for a change of orientation)
and the summand 
$$\frac{1}{1+\zeta}$$
where $\zeta$ is any polynomial in $x,y_1,y_2,z_1,z_2$ by
$$\frac{1}{1+a/\zeta}$$
where
$$a=xy_1y_2z_1z_2v.$$
Unlisted entries $M_{i,j}$ are defined to be $0$.

\vspace{3mm}

\begin{lemma}
\label{l*}
Consider the field $F(Q,Q^\prime)$ where $Q,Q^\prime$ are algebraically
independent over$F$. Let $\phi:F(Q)\r F(Q,Q^\prime)$ be identity on $F$, and
let $\phi(Q)=Q^\prime$. Let $V$ be the intersection of the row space of $M$ with
$\langle U_6,\dots,_{11}\rangle$ where $U_i$ denotes the row vector with 
$1$ in the column corresponding to $u_i$, and $0$'s in the other columns.
Then $dim_{F(Q)}(V)=3$. Additionally, let $w_1,w_2,w_3$ be a basis of $V$
consisting of vectors for which there exist different $i_1,i_2,i_3\in\{6,\dots,11\}$
such that $w_j$ has $i_k$'th coordinate equal to $\delta_{j}^{k}$, $j,k\in\{1,2,3\}$.
Then
\beg{el+}{r=3-rank_{F(Q,Q^\prime)} 
\left(\begin{array}{c}
w_1-\phi(w_1)\\
w_2-\phi(w_2)\\
w_3-\phi(w_3)
\end{array}
\right).
}
\end{lemma}

\vspace{3mm}
\Proof
The matrix $M$ is a submatrix of the matrix of differentials of $\mc_0(\md_{(3,3,0)}$.
More explicitly, we will index things so that $M$ is the $\{1,\dots,9\}\times \{0,\dots, 11\}$
submatrix of the $\{1,\dots 18\}\times \{0,\dots, 17\}$ matrix $N$. Again, we will
denote the rows of $N$ by $v_i$ and columns by $u_j$.
Explicitly, in $B_{(3,3,0)}$, there are additional spanning trees $T_{2,i,j}$ which are
obtained by replacing $e_0$ by $e_2$ in $T_{0,i,j}$. These correspond to
additional columns
$$\begin{array}{lll}
u_{12}=T_{2,2,3}, & u_{13}=T_{2,2,2}, & u_{14}=T_{2,2,1},\\
u_{15}=T_{2,1,1}, & u_{16}=T_{2,1,2}, & u_{17}=T_{2,1,3}.
\end{array}
$$
There are also $9$ additional spanning trees $T^{\prime\prime}_{i,j}$ 
obtained by replacing in $T^{\prime}_{i,j}$ the edge $e_0$ by $e_2$.
Let the row $v_{9+i}$, $i=1,\dots, 9$, be obtained from the row $v_i$
by replacing $T^\prime$ with $T^{\prime\prime}$, thus obtaining $9$ additional rows.
The additional non-zero entries of $N$ are described as follows:
The $(i+9,j+12)$-entry ($i=1,\dots,9$, $j=0,\dots,11$) is obtained from
the $(i,j)$-entry by replacing $A$, $B$ with $A^\prime$, $B^\prime$
where, in the field
$$K^\prime =F(Q^\prime),$$
$$A^\prime=\fracd{1}{1+Q^\prime},\; B^\prime=\fracd{Q^\prime}{Q^\prime +T}.$$
(One has $Q^\prime=f_1g$, $T/Q^\prime=f_2/g$ where $g$ is the face
between $e_0$ and $e_2$.)

Now let $N_1$ resp. $N_2$ be the  $\{1,\dots 9\}\times \{0,\dots, 17\}$
resp. $\{10,\dots,18\}\times \{0,\dots,17\}$ submatrices.
First note that the rank of each of the matrices $N_1$ and $N_2$ is $9$ by the
calculation of the BOS cohomology of $8_{19}$ in \cite{kk} (it is also verified by the
computer-assisted calculation which we will describe below).
This implies that the space $V$ defined in the statement of the Lemma has
$$dim_{F(Q)}(V)=3.$$
Now let $w_1,w_2,w_3$ be a basis as in the statement of the Lemma. By equality
of row and column rank, $r$ is the rank of the $F(Q,Q^\prime)$-space of $6-tuples$
$(\alpha_1,\dots,\alpha_6)\in F(Q,Q^\prime)^6$ such that
$$\alpha_1 w_1+\alpha_2w_2+\alpha_3w_3=
\alpha_4\phi(w_1)+\alpha_5\phi(w_2)+\alpha_6\phi(w_3).$$
Obviously, however, by the assumptions about $w_i$, we must have
$$\alpha_i=\alpha_{3+i},\; i=1,2,3,$$
and the statement follows. 
\qed

\vspace{3mm}
Now to use the Lema, we first construct explicitly a non-zero element $w\in V$ which is of 
the form
$$w=\sum_{i=6}^{11}\alpha_iu_i,\; \alpha_i\in F_0.$$
This will show $r\geq 1$. To construct $w$, let
$$\begin{array}{l}
X=V_1+V_2+V_4+V_5,\\
Y=V_1+V_2+V_7+V_8,\\
Z=V_1+V_3+V_4+V_6,\\
T=V_1+V_3+V_7+V_9
\end{array}
$$
where we denote by $V_i$ the $i$'th row vector of $M$.
One sees immediately from the definition of the row vectors $V_i$
that $X,Y,Z,T$ are linear combinations of the vectors $U_j$, $j=0,\dots,11$
with coefficients in $F_0$.

Now putting
$$\begin{array}{ll}p_1=M_{4,1}+M_{5,1} & q_1=M_{2,4}+M_{5,4}\\
p_2=M_{7,2}+M_{8,2} & q_2=M_{2,4}+M_{8,4}\\
p_3=M_{4,1}+M_{6,1} & q_3=M_{3,5}+M_{6,5}\\
p_4=M_{7,2}+M_{9,2} & q_4=M_{3,5}+M_{9,5}
\end{array}
$$
($p_1$ and $p_3$ are the $u_1$-coordinates of $X,Z$ respectively,
$p_3$ and $p_4$ are the $u_2$-coordinates of $Y,T$ respectively,
$q_1$ and $q_2$ are the $u_4$-coordinates of $X,Y$ respectively,
and $q_3$ and $q_4$ are the $u_5$-coordinates of $Z,T$ respectively;
those are all the non-zero $u_1,u_2,u_4,u_5$ coordinates of $X,Y,Z,T$).

Then one verifies by hand that
$$\frac{p_1p_4}{p_2p_3}=\frac{q_1q_4}{q_2q_3}.
$$
This means that in the vector
$$w=X+\frac{q_2}{q_1}Y +\frac{p_3}{p_1}Z + \frac{p_4}{p_2}\frac{q_2}{q_1}T,$$
the $u_1,u_2,u_4,u_5$ coordinates vanish. One then checks, by hand again,
that the $u_0$ and $u_3$ coordinates vanish as well, thus proving the desired
statement about $w$.

\vspace{3mm}
Proving that $r\leq 1$ is done by using Lemma \ref{l*}. We did this as follows:
Since we are hoping to detect the absence of a relation at a generic point, it is possible
to work at a special point (since a relation absent at a special point cannot occur at
the generic point, using the argument made in detail at the end of \cite{kk}).
Thus, we re-wrote the matrix $M$ over the ring $R=F_0[A,B]$ where 
$F_0=\mathbb{F}_2(t)$, setting
$$\begin{array}{c}x=t,\;y_1=t^2,\;y_2=T,\; z_1=t^4,\\
z_2=t,\; v=t^6,\; a=t^{15}.
\end{array}$$
(The choices of the exponents are arbitrary, with the understanding that too
special choices could create unwanted special relations; a field of rational
functions in a single variable was chosen because computer algebra systems
seem to work much more efficiently in that setting.) We then used Sage to
execute manually the Buchberger algorithm for finding a Gr\"{o}bner basis
of $\langle V_1,\dots,V_9\rangle$ ver the ring $R$, with lexicographic
ordering $u_0>u_1>\dots>u_{11}>A,B$ and degree-lex $A>B$ order in
$A,B$ (the latter of which was chosen because Sage naturally uses that
ordering when working with $F_0[A,B]$). The main reason we worked 
manually is to be able to use heuristics (such as identifying the vector $w$ above)
for speeding up the algorithm. In $70$ easy steps, the Gr\"{o}bner basis elements
we found had leading terms
$$\begin{array}{c}
u_0,\;u_1,\;u_2,\;u_5,\;AB^2u_6,\;A^2u_6,\\
u_7,\; Au_8,\; B^2 u_8,\; A^2B^2u_9.
\end{array}$$
Note that for our purposes, having a Gr\"{o}bner basis is actually 
irrelevant; again, it is merely a tool for performing Gauss elimination
over the fraction field of $R$ which, when done by brute force, would exceed
the computational power of our current implementation of Sage. 
We may then get $w_1,w_2,w_3$ by taking $w$ and our Gr\"{o}bner
basis vectors with leading terms $Au_8$ and $A^2B^2u_9$ and bringing
them to reduced row echelon form, using the substitution
$$A=\fracd{1}{1+Q},\; B=\fracd{Q}{Q+T},\; T=t^10.$$
Again, the choice of $T$ was arbitrary, hoping to avoid a special relation. 
As it turns out, when construction the reduced row echelon form,
we can actually ignore $w$, since we already know it results in a zero
row. 

We used Sage to find by direct computation that $w_2-\phi(w_2)$, $w_3-\phi(w_3)$
are linearly independent (this took several minutes), thus concluding that $r=1$.
This concludes the proof of Proposition \ref{p1}.

\section{Other links with non-trivial BOS cohomologies and branched double covers
with interesting $\widehat{HF}$-homologies}
\label{s2}

\vspace{3mm}

\begin{proposition}
\label{p1s2}
We have
$$rank(H^i(\mc_0(\md_{(3,3,k)})))=\left\{
\begin{array}{ll}
1 &\text{for $i=1,2$}\\
0 &\text{else.}
\end{array}\right.$$
\end{proposition}

\Proof
We proceed by induction on $k$. For $k=0$, this is the statement of Proposition \ref{p1}.
Suppose the statement is true for a given $k$. Consider the long exact sequence for
the cohomology of $\mc_0(\md_{(3,3,k+1)})$ obtained by resolving the edge
$a_{3,k+1}$. Then the $1$-resolution is actually an unlink with $2$ components, and
hence has $0$ BOS cohomology. On the other hand, the $0$-resolution is
$\md_{(3,3,k)}$. Thus, from \rref{ei1} we obtain
$$H^{i}(\mc_0(\md_{(3,3,k+1)}))\cong
H^i(\mc_0(\md_{(3,3,k)})),$$
and the induction step is complete.
\qed

\vspace{3mm}

\begin{corollary}
\label{c1}
We have
$$rank(H^i(\mc_0(\md_{(3,6)})))=\left\{
\begin{array}{ll}
1 &\text{for $i=1,2$}\\
0 &\text{else.}
\end{array}\right.$$
\end{corollary}

\vspace{3mm}

\Proof
Consider the long exact sequence \rref{ei1} form $\md_{(3,3,3)}$
resolving the edge $e_0$. The $0$-resolution is $\md_{(3,6)}$,
the $1$-resolution is an unlink with $2$ components, hence has trivial
BOS cohomology. We conclude that
$$H^{i}(\mc_0(\md_{(3,3,3)}))\cong
H^i(\mc_0(\md_{(3,6)})),$$
and the statement follows from Proposition \ref{p1s2}.
\qed

\vspace{3mm}
Unfortunately, all the examples of links for which we have computed
non-trivial BOS cohomology so far have determinant $0$,
so we cannot use the BOS spectral sequence to make conclusions about $\widehat{HF}$
of their branched double covers. Consider now the black graph $B_k=B(\mathcal{E}_k)$,
$k\geq 2$, depicted in Figure 3. Denote the corresponding link by $L_k$. This is
a knot if $k\geq 3$ is odd and a link with two components when $k\geq 2$ is even.

\begin{figure}
\begin{tikzpicture}[y=0.80pt, x=0.8pt,yscale=-1, inner sep=0pt, outer sep=0pt]
\path[shift={(1113.3049,-645.76085)},draw=black,fill=black,line join=miter,line
  cap=butt,even odd rule,line width=0.800pt]
  (337.1429,266.1122)arc(0.000:180.000:5.535715 and
  4.821)arc(-180.000:0.000:5.535715 and 4.821) -- cycle;
\path[draw=black,line join=miter,line cap=butt,line width=0.800pt]
  (1448.5558,-382.8607) -- (1549.9843,-455.7178);
\path[draw=black,line join=miter,line cap=butt,line width=0.800pt]
  (1337.1272,-454.2892) -- (1441.4129,-381.4321) -- (1322.8415,-261.4321);
\path[draw=black,line join=miter,line cap=butt,line width=0.800pt]
  (1448.5558,-377.1464) -- (1527.1272,-247.1464);
\path[draw=black,dash pattern=on 9.60pt off 9.60pt,line join=miter,line
  cap=butt,miter limit=4.00,line width=0.800pt] (1445.6986,-385.7178) --
  (1444.2700,-490.0035);
\path[shift={(1218.3049,-722.18942)},draw=black,fill=black,line join=miter,line
  cap=butt,even odd rule,line width=0.800pt]
  (337.1429,266.1122)arc(0.000:180.000:5.535715 and
  4.821)arc(-180.000:0.000:5.535715 and 4.821) -- cycle;
\path[shift={(1194.7335,-514.33227)},draw=black,fill=black,line join=miter,line
  cap=butt,even odd rule,line width=0.800pt]
  (337.1429,266.1122)arc(0.000:180.000:5.535715 and
  4.821)arc(-180.000:0.000:5.535715 and 4.821) -- cycle;
\path[shift={(992.59062,-529.33227)},draw=black,fill=black,line join=miter,line
  cap=butt,even odd rule,line width=0.800pt]
  (337.1429,266.1122)arc(0.000:180.000:5.535715 and
  4.821)arc(-180.000:0.000:5.535715 and 4.821) -- cycle;
\path[shift={(1005.4478,-719.33227)},draw=black,fill=black,line join=miter,line
  cap=butt,even odd rule,line width=0.800pt]
  (337.1429,266.1122)arc(0.000:180.000:5.535715 and
  4.821)arc(-180.000:0.000:5.535715 and 4.821) -- cycle;
\path[shift={(1111.8763,-757.18941)},draw=black,fill=black,line join=miter,line
  cap=butt,even odd rule,line width=0.800pt]
  (337.1429,266.1122)arc(0.000:180.000:5.535715 and
  4.821)arc(-180.000:0.000:5.535715 and 4.821) -- cycle;
\path[draw=black,dash pattern=on 9.60pt off 9.60pt,line join=miter,line
  cap=butt,miter limit=4.00,line width=0.800pt] (1444.2700,-492.8607) --
  (1547.1272,-459.2892) -- (1568.3676,-433.4821);
\path[draw=black,dash pattern=on 9.60pt off 9.60pt,line join=miter,line
  cap=butt,miter limit=4.00,line width=0.800pt] (1529.2700,-248.5749) --
  (1566.4129,-284.2892);
\path[draw=black,dash pattern=on 6.40pt off 6.40pt,line join=miter,line
  cap=butt,miter limit=4.00,line width=0.800pt] (1338.5558,-455.7178) --
  (1441.4129,-490.7178);
\path[draw=black,dash pattern=on 9.60pt off 9.60pt,line join=miter,line
  cap=butt,miter limit=4.00,line width=0.800pt] (1333.5558,-450.0035) --
  (1297.1272,-410.7178);
\path[draw=black,dash pattern=on 9.60pt off 9.60pt,line join=miter,line
  cap=butt,miter limit=4.00,line width=0.800pt] (1321.4129,-264.2892) --
  (1296.4129,-313.5749);
\path[draw=black,dash pattern=on 9.60pt off 9.60pt,line join=miter,line
  cap=butt,miter limit=4.00,line width=0.800pt] (1297.8415,-410.0035) --
  (1295.6986,-313.5750);
\path[shift={(967.56704,-677.41094)},draw=black,fill=black,line join=miter,line
  cap=butt,even odd rule,line width=0.800pt]
  (337.1429,266.1122)arc(0.000:180.000:5.535715 and
  4.821)arc(-180.000:0.000:5.535715 and 4.821) -- cycle;
\path[shift={(964.53658,-580.4363)},draw=black,fill=black,line join=miter,line
  cap=butt,even odd rule,line width=0.800pt]
  (337.1429,266.1122)arc(0.000:180.000:5.535715 and
  4.821)arc(-180.000:0.000:5.535715 and 4.821) -- cycle;
\path[shift={(1054.4402,-504.16978)},draw=black,fill=black,line join=miter,line
  cap=butt,even odd rule,line width=0.800pt]
  (337.1429,266.1122)arc(0.000:180.000:5.535715 and
  4.821)arc(-180.000:0.000:5.535715 and 4.821) -- cycle;
\path[shift={(1131.2117,-500.12917)},draw=black,fill=black,line join=miter,line
  cap=butt,even odd rule,line width=0.800pt]
  (337.1429,266.1122)arc(0.000:180.000:5.535715 and
  4.821)arc(-180.000:0.000:5.535715 and 4.821) -- cycle;
\path[draw=black,dash pattern=on 9.60pt off 9.60pt,line join=miter,line
  cap=butt,miter limit=4.00,line width=0.800pt] (1325.8252,-260.8361) --
  (1384.9191,-238.1077) -- (1461.1857,-233.5620) -- (1526.8456,-246.1889);
\path[shift={(1235.9826,-698.95591)},draw=black,fill=black,line join=miter,line
  cap=butt,even odd rule,line width=0.800pt]
  (337.1429,266.1122)arc(0.000:180.000:5.535715 and
  4.821)arc(-180.000:0.000:5.535715 and 4.821) -- cycle;
\path[shift={(1235.4775,-552.48379)},draw=black,fill=black,line join=miter,line
  cap=butt,even odd rule,line width=0.800pt]
  (337.1429,266.1122)arc(0.000:180.000:5.535715 and
  4.821)arc(-180.000:0.000:5.535715 and 4.821) -- cycle;
\path[draw=black,dash pattern=on 6.40pt off 6.40pt,line join=miter,line
  cap=butt,miter limit=4.00,line width=0.800pt] (1567.7567,-287.1001) --
  (1571.2923,-306.2930);
\path[draw=black,dash pattern=on 0.68pt off 0.68pt,line join=miter,line
  cap=butt,miter limit=4.00,line width=0.682pt] (1580.3837,-375.7480) --
  (1571.7974,-311.8488);
\path[shift={(1250.1247,-662.08534)},draw=black,fill=black,line join=miter,line
  cap=butt,even odd rule,line width=0.800pt]
  (337.1429,266.1122)arc(0.000:180.000:5.535715 and
  4.821)arc(-180.000:0.000:5.535715 and 4.821) -- cycle;
\path[draw=black,dash pattern=on 9.60pt off 9.60pt,line join=miter,line
  cap=butt,miter limit=4.00,line width=0.800pt] (1582.9090,-392.1560) --
  (1581.3938,-377.0037);
\path[draw=black,dash pattern=on 9.60pt off 9.60pt,line join=miter,line
  cap=butt,miter limit=4.00,line width=0.800pt] (1568.7669,-429.5316) --
  (1582.9090,-395.6915);
\path[xscale=0.947,yscale=1.056,fill=black] (1573.8591,-464.22076) node[above
  right] (text4512) {$a$     };
\path[xscale=0.798,yscale=1.253,fill=black] (1952.7096,-356.22076) node[above
  right] (text4516) {$e_1$     };
\path[xscale=0.919,yscale=1.089,fill=black] (1719.9667,-379.21136) node[above
  right] (text4520) {$e_2$     };
\path[xscale=0.665,yscale=1.504,fill=black] (2317.9697,-161.37625) node[above
  right] (text4524) {$e_k$   };

\end{tikzpicture}
\caption{}
\end{figure}

\vspace{3mm}

\begin{theorem}
\label{t1}For $k\geq 2$,
$$rank H^i(\mc_0(\mathcal{E}_k))=
\left\{
\begin{array}{ll}
1 & \text{for $i=3/2,5/2$}\\
k-2 & \text{for $i=7/2$}\\
0 & \text{else.}
\end{array}
\right.
$$
\end{theorem} 

\vspace{3mm}
\Proof
Resolve the projection $\mathcal{E}_k$ at the edge $a$.
The $0$-resolution $\mathcal{E}^{0}_{k}$ is $\mathcal{D}_{(3,3,k+1)}$ after
undoing a single R2 move, thus,
$$H^i(\mc_0(\md_{(3,3,k+1)}))\cong H^{i+1/2}(\mc_0(\mathcal{E}^{0}_{k})),$$
i.e.
$$rank H^i(\mc_0(\mathcal{E}^{0}_{k}))=
\left\{
\begin{array}{ll}
1 & \text{for $i=3/2,5/2$}\\
0 & \text{else.}
\end{array}
\right.
$$
On the other hand, the $1$-resolution can be processed as follows: An R3 move
combined with undoing a positive (non-height changing) R1 move gives a move
shown in Figure 4.

\begin{figure}
\begin{tikzpicture}[y=0.80pt, x=0.8pt,yscale=-1, inner sep=0pt, outer sep=0pt]
\path[draw=black,line join=miter,line cap=butt,line width=0.800pt]
  (-1644.5283,763.4585) -- (-1570.7872,832.1489);
\path[shift={(-1978.8859,494.42505)},draw=black,fill=black,line join=miter,line
  cap=butt,even odd rule,line width=0.800pt]
  (337.1429,266.1122)arc(0.000:180.000:5.535715 and
  4.821)arc(-180.000:0.000:5.535715 and 4.821) -- cycle;
\path[shift={(-1902.1143,565.13573)},draw=black,fill=black,line join=miter,line
  cap=butt,even odd rule,line width=0.800pt]
  (337.1429,266.1122)arc(0.000:180.000:5.535715 and
  4.821)arc(-180.000:0.000:5.535715 and 4.821) -- cycle;
\path[shift={(-1902.1143,565.13573)},draw=black,fill=black,line join=miter,line
  cap=butt,even odd rule,line width=0.800pt]
  (337.1429,266.1122)arc(0.000:180.000:5.535715 and
  4.821)arc(-180.000:0.000:5.535715 and 4.821) -- cycle;
\path[shift={(-1879.8909,475.23215)},draw=black,fill=black,line join=miter,line
  cap=butt,even odd rule,line width=0.800pt]
  (337.1429,266.1122)arc(0.000:180.000:5.535715 and
  4.821)arc(-180.000:0.000:5.535715 and 4.821) -- cycle;
\path[draw=black,dash pattern=on 9.60pt off 9.60pt,line join=miter,line
  cap=butt,miter limit=4.00,line width=0.800pt] (-1650.5893,759.4179) --
  (-1551.5943,741.2352) -- (-1571.7974,829.1185);
\path[draw=black,line join=miter,line cap=butt,line width=0.800pt]
  (-1646.5486,761.4382) -- (-1706.1477,759.4179);
\path[draw=black,dash pattern=on 9.60pt off 9.60pt,line join=miter,line
  cap=butt,miter limit=4.00,line width=0.800pt] (-1686.9548,782.6514) --
  (-1650.5893,760.4281) -- (-1666.7517,797.8037);
\path[draw=black,dash pattern=on 9.60pt off 9.60pt,line join=miter,line
  cap=butt,miter limit=4.00,line width=0.800pt] (-1614.2238,853.3621) --
  (-1571.7974,833.1591) -- (-1588.9700,863.4637);
\path[draw=black,line join=miter,line cap=butt,line width=0.800pt]
  (-1567.7567,834.1692) -- (-1569.7771,875.5855);
\path[draw=black,dash pattern=on 6.40pt off 6.40pt,line join=miter,line
  cap=butt,miter limit=4.00,line width=0.800pt] (-1705.1375,759.4179) --
  (-1771.8075,761.4382);
\path[draw=black,line join=miter,line cap=butt,line width=0.800pt]
  (-1695.0360,786.6920) -- (-1758.6756,825.0779);
\path[draw=black,line join=miter,line cap=butt,line width=0.800pt]
  (-1670.7924,804.8748) -- (-1699.0766,853.3621);
\path[draw=black,line join=miter,line cap=butt,line width=0.800pt]
  (-1626.3456,858.4129) -- (-1673.8228,883.6667);
\path[draw=black,line join=miter,line cap=butt,line width=0.800pt]
  (-1595.0309,872.5550) -- (-1626.3456,913.9713);
\path[draw=black,dash pattern=on 9.60pt off 9.60pt,line join=miter,line
  cap=butt,miter limit=4.00,line width=0.800pt] (-1570.7872,883.6667) --
  (-1581.8989,941.2454);
\path[shift={(-1630.3832,500.48597)},draw=black,fill=black,line join=miter,line
  cap=butt,even odd rule,line width=0.800pt]
  (337.1429,266.1122)arc(0.000:180.000:5.535715 and
  4.821)arc(-180.000:0.000:5.535715 and 4.821) -- cycle;
\path[shift={(-1553.6116,571.19665)},draw=black,fill=black,line join=miter,line
  cap=butt,even odd rule,line width=0.800pt]
  (337.1429,266.1122)arc(0.000:180.000:5.535715 and
  4.821)arc(-180.000:0.000:5.535715 and 4.821) -- cycle;
\path[shift={(-1553.6116,571.19665)},draw=black,fill=black,line join=miter,line
  cap=butt,even odd rule,line width=0.800pt]
  (337.1429,266.1122)arc(0.000:180.000:5.535715 and
  4.821)arc(-180.000:0.000:5.535715 and 4.821) -- cycle;
\path[shift={(-1531.3882,481.29307)},draw=black,fill=black,line join=miter,line
  cap=butt,even odd rule,line width=0.800pt]
  (337.1429,266.1122)arc(0.000:180.000:5.535715 and
  4.821)arc(-180.000:0.000:5.535715 and 4.821) -- cycle;
\path[draw=black,line join=miter,line cap=butt,miter limit=4.00,line
  width=0.800pt] (-1302.0866,765.4789) -- (-1203.0917,747.2961) --
  (-1223.2947,835.1794);
\path[draw=black,line join=miter,line cap=butt,line width=0.800pt]
  (-1298.0460,767.4992) -- (-1357.6450,765.4789);
\path[draw=black,dash pattern=on 9.60pt off 9.60pt,line join=miter,line
  cap=butt,miter limit=4.00,line width=0.800pt] (-1338.4521,788.7124) --
  (-1302.0866,766.4890) -- (-1318.2490,803.8646);
\path[draw=black,dash pattern=on 9.60pt off 9.60pt,line join=miter,line
  cap=butt,miter limit=4.00,line width=0.800pt] (-1265.7211,859.4230) --
  (-1223.2947,839.2200) -- (-1240.4673,869.5246);
\path[draw=black,line join=miter,line cap=butt,line width=0.800pt]
  (-1219.2541,840.2301) -- (-1221.2744,881.6464);
\path[draw=black,dash pattern=on 6.40pt off 6.40pt,line join=miter,line
  cap=butt,miter limit=4.00,line width=0.800pt] (-1356.6348,765.4788) --
  (-1423.3048,767.4992);
\path[draw=black,line join=miter,line cap=butt,line width=0.800pt]
  (-1346.5333,792.7530) -- (-1410.1729,831.1388);
\path[draw=black,line join=miter,line cap=butt,line width=0.800pt]
  (-1322.2897,810.9357) -- (-1350.5739,859.4230);
\path[draw=black,line join=miter,line cap=butt,line width=0.800pt]
  (-1277.8430,864.4738) -- (-1325.3201,889.7276);
\path[draw=black,line join=miter,line cap=butt,line width=0.800pt]
  (-1246.5282,878.6159) -- (-1277.8430,920.0322);
\path[draw=black,dash pattern=on 9.60pt off 9.60pt,line join=miter,line
  cap=butt,miter limit=4.00,line width=0.800pt] (-1222.2846,889.7276) --
  (-1233.3963,947.3063);
\path[draw=black,line join=miter,line cap=butt,miter limit=4.00,line
  width=0.798pt] (-1433.9065,855.3824) -- (-1528.3608,853.3621);
\path[draw=black,line join=miter,line cap=butt,line width=0.800pt]
  (-1470.7821,822.0474) -- (-1432.3963,854.3723) -- (-1466.7415,880.6363);

\end{tikzpicture}
\caption{}
\end{figure}

Undoing three R2 moves and three positive R1 moves, we obtain the black graph
shown in Figure 5 and, after undoing two R2 moves, we obtain
a cycle of $k-2$ height $0$ edges when $k>2$, and an unlink of
two components when $k=2$.

\begin{figure}
\begin{tikzpicture}[y=0.80pt, x=0.8pt,yscale=-1, inner sep=0pt, outer sep=0pt]
\path[draw=black,line join=miter,line cap=butt,line width=0.800pt]
  (-165.4365,369.1401) -- (-247.2588,447.9320) -- (-171.4974,503.4904);
\path[shift={(-499.53649,240.89057)},draw=black,fill=black,line join=miter,line
  cap=butt,even odd rule,line width=0.800pt]
  (337.1429,266.1122)arc(0.000:180.000:5.535715 and
  4.821)arc(-180.000:0.000:5.535715 and 4.821) -- cycle;
\path[shift={(-499.53649,240.89057)},draw=black,fill=black,line join=miter,line
  cap=butt,even odd rule,line width=0.800pt]
  (337.1429,266.1122)arc(0.000:180.000:5.535715 and
  4.821)arc(-180.000:0.000:5.535715 and 4.821) -- cycle;
\path[shift={(-494.53649,102.31914)},draw=black,fill=black,line join=miter,line
  cap=butt,even odd rule,line width=0.800pt]
  (337.1429,266.1122)arc(0.000:180.000:5.535715 and
  4.821)arc(-180.000:0.000:5.535715 and 4.821) -- cycle;
\path[draw=black,dash pattern=on 9.60pt off 9.60pt,line join=miter,line
  cap=butt,miter limit=4.00,line width=0.800pt] (-161.8365,369.1332) --
  (-102.5508,397.7046);
\path[draw=black,dash pattern=on 9.60pt off 9.60pt,line join=miter,line
  cap=butt,miter limit=4.00,line width=0.800pt] (-166.8365,506.2760) --
  (-114.6937,489.1332);
\path[shift={(-436.67939,131.60485)},draw=black,fill=black,line join=miter,line
  cap=butt,even odd rule,line width=0.800pt]
  (337.1429,266.1122)arc(0.000:180.000:5.535715 and
  4.821)arc(-180.000:0.000:5.535715 and 4.821) -- cycle;
  
  \path[shift={(-578.67939,180.60485)},draw=black,fill=black,line join=miter,line
  cap=butt,even odd rule,line width=0.800pt]
  (337.1429,266.1122)arc(0.000:180.000:5.535715 and
  4.821)arc(-180.000:0.000:5.535715 and 4.821) -- cycle;
\path[shift={(-443.10799,221.60485)},draw=black,fill=black,line join=miter,line
  cap=butt,even odd rule,line width=0.800pt]
  (337.1429,266.1122)arc(0.000:180.000:5.535715 and
  4.821)arc(-180.000:0.000:5.535715 and 4.821) -- cycle;
\path[draw=black,dash pattern=on 0.53pt off 3.17pt,line join=miter,line
  cap=butt,miter limit=4.00,line width=0.529pt] (-105.4080,401.3204) --
  (-112.5508,486.2760);
\path[xscale=0.838,yscale=1.193,fill=black] (-177.30659,312.16937) node[above
  right] (text5725) {$e_1$     };
\path[xscale=0.676,yscale=1.478,fill=black] (-227.73466,357.13852) node[above
  right] (text5729) {$e_k$     };

\end{tikzpicture}
\caption{}
\end{figure}

We have then
$$rank H^i(\mc_0(\mathcal{E}^{1}_{k}))=
\left\{
\begin{array}{ll}
k-2 & \text{for $i=3$}\\
0 & \text{else.}
\end{array}
\right.
$$
For $k=2$, we are therefore done. For $k>2$, we are done if we can show the
following statement. 
\qed

\begin{lemma}
\label{l2}
The connecting map
$$\diagram
H^i(\mc_0(\mathcal{E}^{0}_{k}))\rto^{\delta} &
H^{i+1/2}(\mc_0(\mathcal{E}^{1}_{k}))
\enddiagram$$
of \rref{ei1} is $0$ for all $i$.
\end{lemma}

\Proof
For $k=3$, $\mathcal{E}_{3}$ is actually a projection of the mirror of
$T(7,3)$, whose 
$\widehat{HF}$ has rank $3$ by \cite{ost3}. Therefore, we must have
$\delta=0$ by the Baldwin-Ozsv\'{a}th-Szab\'{o} spectral sequence.

Now our proof will be by induction on $k\geq 3$. Consider, for $k>3$,
the $0$-resolutions $\mathcal{E}^{(1)}_{k}$, $\mathcal{E}^{(2)}_{k}$
of $\mathcal{E}_{k}$ at $e_{k-1}$, $e_{k}$. Then every spanning tree
of the black graph of $\mathcal{E}_{k}$ gives rise to a spanning tree
of  $\mathcal{E}^{(i)}_{k}$ for $i=1$ or $i=2$, and the spanning trees which
give rise to both give rise to spanning trees of the $0$-resolution $\mathcal{E}^{(12)}_{k}$
at both $e_{k-1}$, $e_{k}$.

We have, then, a ``Mayer-Vietoris exact sequence''
\beg{emv+}{\diagram
0\rto & mc_0(\mathcal{E}_k)\rto^(.35)\iota &
\mc_0(\mathcal{E}^{(1)}_{k})\oplus
\mc_0(\mathcal{E}^{(2)}_{k})\rto &
\mc_0(\mathcal{E}^{(12)}_{k})\rto & 0.
\enddiagram
}
One has, of course, 
$$\mc_0(\mathcal{E}^{(i)}_{k})\cong \mc_0(\mathcal{E}_{k-1}),$$
$$\mc_0(\mathcal{E}^{(12)}_{k})\cong \mc_0(\mathcal{E}_{k-2}).$$
Moreover, the maps \rref{emv+} induce maps of the long exact sequences
corresponding to resolution at the edge $a$. In particular, we obtain
a commutative square
\beg{emv*}{\diagram\protect
H^i(\mc_0(\mathcal{E}^{0}_{k}))\dto_{\iota_*}\rto^\delta &
H^{i+1/2}(\mc_0(\mathcal{E}^{1}_{k}))\dto^{\iota_*}\\
{\protect\begin{array}{c}
H^i(\mc_0(\mathcal{E}^{(1)0}_{k}))\\
\oplus\\
H^i(\mc_0(\mathcal{E}^{(2)0}_{k}))
\end{array}}
\rto^{\delta\oplus\delta} &
{\protect\begin{array}{c}
H^{i+1/2}(\mc_0(\mathcal{E}^{(1)1}_{k}))\\
\oplus\\
H^{i+1/2}(\mc_0(\mathcal{E}^{(2)1}_{k})).
\end{array}}
\enddiagram
}
By the induction hypothesis, the bottom row satisfies $\delta\oplus\delta=0$,
while the left column of \rref{emv*} is injective by our computation of the
$0$-resolutions (the two components omit the $(k-3)$'rd and $(k-2)$'nd
summands of $F^{\oplus k}$ where $F$ is the ground field, respectively).
Since the left column of \rref{emv*} is injective, the top row then
satisfies $\delta=0$, as claimed. This concludes the proof of the Lemma and
hence the Theorem.
\qed

\vspace{3mm}

\begin{theorem}
\label{t4}
For $k\geq 3$,  we have
$$rank \widehat{HF}(\Sigma(L_k))=k,$$
while
$$det(L_k)=k-2.$$
\end{theorem}

\Proof
The computation of the determinant follows from Theorem \ref{t1} (since the
determinant is, up to sign, the trace of BOS cohomology).
Since $det(L_k)\neq 0$, the BOS spectral sequence then applies, with the
$E_3$-term given by Theorem \ref{t1}. By sparsity, no differential is possible,
and hence the spectral sequence collapses to $E_3$ in this case.
\qed

\vspace{3mm}
\noindent
{\bf Comment:} Note that while Theorems \ref{t1}, \ref{t4} do not
provide examples of non-collapse of the BOS spectral sequence,
they exhibit interesting behavior in the sense
of an ``extension'': The BOS cohomology of $L_k$ has
non-trivial elements in degrees $3/2$ and $7/2$, which are congruent modulo $2$.

\vspace{3mm}
\begin{proposition}
\label{thyp}
For all but finitely many values of $k>3$, the link $L_k$ (knot when $k$ is odd)
is hyperbolic
\end{proposition}

\Proof
The moves converting $\mathcal{E}_3$ to the mirror of the standard knot projection
of $T(7,3)$ can be made in such
a way that the crossing $x$ corresponding to the edge $e_2$ in Figure 3 is
not involved in any Reidemeister move. Form a link
$M_3$ by adding an unknotted link component $\ell$ to $L_3$ encircling the
crossing $x$. Using SnapPea, the link $M_3$ is hyperbolic
with volume $6.551743287888$.
Now $L_k$ for $k>3$ can be obtained from $M_3$ by
performing hyperbolic Dehn filling on the link component $\ell$.
Because of this, all but finitely many of the links $L_k$ are hyperbolic by
Thurston's theorem \cite{thur} (see also \cite{a}, Section 3).
\qed

\vspace{3mm}
\noindent
{\bf Comment:} The only example of $k>3$ we know for which 
$L_k$ is not hyperbolic is $k=5$. The knot $L_5$ is actually the mirror
image of $T(8,3)$. The Jones polynomial of the mirror of $L_k$, $k\geq 3$,
is
$$t^{(k+9)/2}(1+t^2-t^9\frac{1+t^{k-4}}{1+t}).$$


\vspace{10mm}

\begin{thebibliography}{99}

\bibitem{a} C.Adams: Hyperbolic knots, {\em in: Handbook of
knot theory}, 1-18, Elsevier, Amsterdam, 2005

%\bibitem{ac} C. Adams: Toroidally alternating knots and links, {\em Topology} 33 (1994) 353-369

\bibitem{os} J.Baldwin, P.Ozsv\'{a}th, Z.Szab\'{o}: Heegaard Floer homology of
double-covers, Kauffman states, and Novikov rings, to appear


\bibitem{kk} D.Kriz, I.Kriz: Baldwin-Ozsv\'{a}th-Szab\'{o} cohomology
is a link invariant, arXiv: 1109.0064, to appear in {\em Advances in Math.}

\bibitem{ls1} P.Lisca, A.Stipsicz: Ozsv\'{a}th-Szab\'{o} invariants and
tight contact 3-manifolds III, {\em J. Sympl. Geom.} 5 (2007) 357-384

\bibitem{ls2} P.Lisca, A.Stipsicz: On the existence of tight contact structures
on Seifert fibered 3-manifolds, {\em Duke Math. J.} 148 (2009) 175-209

\bibitem{osp} P.Ozsv\'{a}th, Z.Szab\'{o}: On the Floer homology
of plumbed three-manifolds, {\em Geom. Topol. 7} (2003) 185-224

\bibitem{os2} P.Ozsv\'{a}th, Z.Szab\'{o}: On the Heegaard Floer homology 
of branched double-covers, {\em Adv. Math.}  194  (2005),  no. 1, 1-33


\bibitem{ost3} P.Ozsv\'{a}th, Z.Szab\'{o}: Holomorphic disks, link invariants 
and the multi-variable Alexander polynomial.  {\em Algebr. Geom. Topol.}  
8  (2008),  no. 2, 615-692

\bibitem{thur} W. Thurston: {\em The geometry and topology of 
$3$-manifolds}, lecture notes, Princeton University, 1978


\end{thebibliography}
\end{document}